# Γ-limit of the cut functional on dense graph sequences.

ANDREA BRAIDES*, PAOLO CERMELLI†, AND SIMONE DOVETTA†‡

**Abstract.** A sequence of graphs with diverging number of nodes is a dense graph sequence if the number of edges grows approximately as for complete graphs. To each such sequence a function, called graphon, can be associated, which contains information about the asymptotic behavior of the sequence. Here we show that the problem of subdividing a large graph in communities with a minimal amount of cuts can be approached in terms of graphons and the Γ-limit of the cut functional, and discuss the resulting variational principles on some examples. Since the limit cut functional is naturally defined on Young measures, in many instances the partition problem can be expressed in terms of the probability that a node belongs to one of the communities. Our approach can be used to obtain insights into the bisection problem for large graphs, which is known to be NP-complete.

**1. Introduction.** Large graphs are ubiquitous in applications to computer networks, power grids, social networks, systems biology, combinatorial optimization, statistics, and many other fields, but their computational treatment is a major challenge, even though computing power available to researchers is steadily increasing. An example is the *bisection problem*, i.e., the problem of subdividing a graph in two parts with the same number of nodes and with minimal number of edges connecting the two parts: it is a basic operation in handling graphs, for instance in divide-and-conquer algorithms, but it is known to be NP-complete [1].

One possibility to study large graphs is to characterize the asymptotic properties of sequences of graphs that grow by some iterative rule, a popular example of this technique being the preferential attachment model for the web graph proposed in [2]. This approach requires to control the properties of each finite graph in the sequence, and does not rely on the notion of a limit object of a sequence of graphs.

An alternative technique to deal with sequences of graphs has been developed recently: the idea is to associate to each such sequence a limit object, called *graphon* ([3], [4], [5], [6], [7], [8], [9], [10]). This is a function on $[0,1]^2$, which in some sense represents the limit of the adjacency matrices of the graphs. The graphon inherits some average topological property of the sequence of graphs that generates it, such as the number of copies of subgraphs. It turns out that each graph sequence admits a non-trivial limit graphon. Moreover, if graphs are suitably identified with some graphons through their adjacency matrices, convergence of graph sequences can be reformulated in terms of the convergence of graphons in a suitable norm, the so-called *cut norm*, defined below.

In this work we study a generalization of the bisection problem for a dense graph sequence, i.e., how to partition a large graph in a finite number of fixed-size subsets with a minimal number of cuts: we show that the problem can be formulated in terms of graphons, and that the limit problem always has solution. To illustrate our approach, we focus below on the bisection problem only: indeed, it is well known that the bisection problem for each graph in the sequence can be restated as a minimum problem for the so-called *cut functional*, which counts the number of edges connecting two regions of the graph. In its simplest form, the cut functional can be written as:

$$\tilde{F}_n(\tilde{u}) = \frac{1}{n^2} \sum_{i,j \in V(G_n)} A^n_{ij} |\tilde{u}(i) - \tilde{u}(j)|^2, \tag{1}$$

where, for every $n \in \mathbb{N}$, $G_n$ is a graph with nodes $V(G_n)$ and edges $E(G_n)$ respectively. Here $A^n_{ij}$ denotes the *adjacency matrix* of $G_n$, i.e., the square matrix such that $A^n_{ij} = 1$ if $ij \in E(G_n)$ and $A^n_{ij} = 0$ otherwise, while $\tilde{u} : V(G_n) \to \{-1, +1\}$ is a *spin function* on $G_n$. We assume that the graph sequence is a *dense graph sequence*, i.e., the number of nodes and edges of $G_n$ grow as $n$ and $n^2$, respectively.

---

* Dipartimento di Matematica, Università di Roma "Tor Vergata", Roma, Italy (braides@mat.uniroma2.it).

† Dipartimento di Matematica, Università di Torino, Via Carlo Alberto, 10, 10123 Torino, Italy (paolo.cermelli@unito.it).

‡ Dipartimento di Scienze Matematiche "G.L. Lagrange", Politecnico di Torino, Corso Duca degli Abruzzi, 24, 10129 Torino, Italy (simone.dovetta@unito.it).

Specifically, we will consider the limit as $n \to +\infty$ of the minimization

$$\min_{\tilde{u} \in \tilde{E}^n} \tilde{F}_n(\tilde{u}) \tag{2}$$

with

$$\tilde{E}^n := \Big\{ \tilde{u} : V(G_n) \to \{-1, +1\} : \sum_{i \in V(G_n)} \tilde{u}(i) = 0 \Big\}.$$

These minimum problems have a natural interpretation in terms of the bisection problem for the graphs $G_n$. Indeed, letting

$$S_n = \{i \in V(G_n) : \tilde{u}(i) = +1\}, \qquad S^c_n = V(G_n) \setminus S_n,$$

then the only non-vanishing terms in (1) are those that involve pairs of nodes $i, j$ such that either $i \in S_n$ and $j \in S^c_n$ or $j \in S_n$ and $i \in S^c_n$, so that, granted the symmetry of the adjacency matrix, we can write

$$\tilde{F}_n(\tilde{u}) = \frac{8}{n^2} \sum_{i \in S_n, j \in S^c_n} A^n_{ij} = \frac{8}{n^2} e_{G_n}(S_n, S^c_n), \tag{3}$$

where $e_{G_n}(S_n, S^c_n)$ is the number of edges that connect $S_n$ to $S^c_n$. Hence, the minimization of the functional (1) over spin functions in $\tilde{E}^n$ is equivalent to finding the subdivision of $V(G_n)$ into equal parts that minimizes the size of the cut between the two communities $S_n$ and $S^c_n$.

To study the asymptotic behaviour of the minimization problems (2) we compute the Γ-limit of the sequence of cut functionals (1). Indeed, Γ-convergence guarantees that the minima of the cut functionals converge to the minima of their limit (see for instance [11], [12], [13], [14], [15], [16]). Our problem belongs to the class of discrete-to-continuum problems, for which Γ-convergence is a standard tool, and indeed Γ-limits of functionals on graph sequences have been calculated ([17], [18], [19]). However, in the existing literature, most graphs are either lattices or are embedded in an Euclidean space, and the kernels of the functionals (the equivalent of $A^n_{ij}$ in (1)) are required to decay with the Euclidean distance between the nodes: the topology of the graph, and the coupling between the nodes is inherited by the metric embedding. On the other hand, in general, the topology of an abstract graph is not dictated by any embedding, and therefore all information regarding connectivity and the interactions between the nodes is contained in the adjacency matrix, which makes the problem truly non local. The advantage of working in terms of graphons is just that no embedding of the nodes has to be performed.

Our main result is the representation, in terms of Young measures, of the Γ-limit of a sequence of cut functionals of quite general form, that allows to study the multi-partition problem in the dense-graph limit. Referring to the example above, the statement is that the Γ-limit of functionals (1) admits the representation

$$(4) \qquad I(\nu) = \int_{[0,1]^2} W(x,y) \Big( \int_{\mathbb{R}^2} |\lambda - \mu|^2 d\nu_x(\lambda) d\nu_y(\mu) \Big) dx dy,$$

where the kernel $W(x,y)$ is the *limit graphon* associated to the sequence of graphs $\{G_n\}_{n \in \mathbb{N}}$, and $\nu_x$ is a *spin Young measure* on $[0,1]$, i.e., a probability measure with support on $\{-1,1\}$, parameterized by $x \in [0,1]$. The interval $[0,1]$ plays here the role of the set of nodes of the graph in the limit. The proof itself actually turns out to be very simple, once we characterize the convergence of graphons in such a way that functionals of the form (4) are continuously convergent if the corresponding graphons converge.

The fact that the bisection problem has a solution in terms of a probability measure reflects the fact that, when the graph has a large number of connections, the distinction between the communities can be blurry, and what makes sense is just to compute the probability that a node belongs to one of the two communities. We discuss how to use the representation (4) of the Γ-limit functional to actually compute bisections of very large graphs that minimize the cut size.

We remark that the integral representations of the limit cut functional that we use here are akin to nonlocal functionals recently studied in the context of generalized theories of continua ([20], [21], [22], [23], [24]).

The paper is organized as follows. Section 2 provides an overview of the underlying theory of converging graph sequences, while Section 3 collects some basic properties of Young measures. Section 4 is devoted to the proof of our main results. Finally, in Section 5, we investigate the minimization problem for Γ−limit functional on several explicit examples of dense graph sequences.

**2. Converging graph sequences.** In this section we briefly revise the theory of graphons and convergence of graph sequences developed by Borgs, Lovasz *et al.* [6, 7, 9, 10, 25].

In what follows $G$ denotes a simple graph, without loops and multiple edges, with $V(G) = \{1, \ldots, n\}$ the set of vertices (nodes) and $E(G)$ the set of edges; we assume that the graph is undirected, i.e., we identify the edges $ij$ and $ji$. We also denote by $A = (A_{ij})$ the $n \times n$ *adjacency matrix* of $G$, such that $A_{ij} = 1$ if $ij \in E(G)$ and $A_{ij} = 0$ otherwise. Since $G$ is undirected, $A$ is symmetric.

For $F$ and $G$ graphs, we say that a map $\phi : V(F) \to V(G)$ is a *homomorphism* if it is adjacency preserving (i.e., if $ij \in E(F)$ then $\phi(i)\phi(j) \in E(G)$), and denote by $\hom(F, G)$ the number of such maps. We define the *homomorphism densities* of $F$ into $G$ as

$$(5) \qquad t(F, G) = \frac{\hom(F, G)}{|V(G)|^{|V(F)|}},$$

where $|\cdot|$ denotes here the cardinality of a finite set[1]. Note that (5) can be interpreted as the normalized number of copies of $F$ that can be found in $G$.

---

[1] In what follows we use the notation $|\cdot|$ for the absolute value of a real number, the cardinality of a finite set and the Lebesgue measure of a measurable subset of $\mathbb{R}^n$, the meaning being clear from the context.

DEFINITION 1. *A sequence of graphs $\{G_n\}_{n\in\mathbb{N}}$ is said to be* left convergent *if the sequences $t(F, G_n)$ converge for $n \to +\infty$, for every simple graph $F$.*

Left convergence ensures that the local structure of the graphs in the sequence asymptotically stabilizes in some sense.

It turns out that there is an object that describes left-converging graph sequences. A *graphon* is a bounded measurable function $W : [0, 1]^2 \to \mathbb{R}$ that is symmetric, i.e., $W(x, y) = W(y, x)$ for every $(x, y) \in [0, 1]^2$ (cf. [5], [25]). We denote by $\mathcal{W}$ the set of all such functions, and by $\mathcal{W}_0$ the set of graphons with values in $[0, 1]$. Heuristically, a graphon is a weighted graph with $[0, 1]$ as the set of nodes, with adjacency/weight matrix given by $W(x, y)$. Standard definitions for graphs can be adapted to graphons: for instance, we say that the *degree* of $x \in [0, 1]$, as a node of $W$, is $\deg_W(x) = \int_{[0,1]} W(x, y) dy$.

It is possible to associate to every graph $G$ a family of graphons. If $|V(G)| = n$, consider any labeling of the nodes, divide the interval $[0, 1]$ in $n$ intervals $I_1^n = [0, \frac{1}{n}]$, $I_2^n = (\frac{1}{n}, \frac{2}{n}]$, ......, $I_n^n = (\frac{(n-1)}{n}, 1]$, and define the piecewise-constant function $W_G$ on $[0, 1]^2$ by

(6) $$W_G(x, y) = A_{ij} \quad \text{if} \quad (x, y) \in I_i^n \times I_j^n.$$

Therefore, $W_G$ is a functional representation of the adjacency matrix of $G$, once its nodes have been labeled someway. Different labelings yield different graphons associated to the same graph.

We define the *cut norm* of $W \in \mathcal{W}$ as

(7) $$\|W\|_\square = \sup_{S,T \subseteq [0,1]} \left| \int_{S \times T} W(x, y) dx dy \right|,$$

with the supremum taken over all measurable subsets $S, T$ of $[0, 1]$. The introduction of the cut norm goes back to Frieze and Kannan [26], but the recognition of its fundamental role in the theory of converging graph sequences is due to Borgs, Chayes, Lovasz, Sos and Vesztergombi (for instance in [6]).

The cut norm can be equivalently defined in any of the following ways ([6]):

$$\|W\|_\square = \sup_{S \subseteq [0,1]} \left| \int_{S \times S^c} W(x, y) dx dy \right|$$
$$= \sup_{S,T \subseteq [0,1], S \cap T = \emptyset} \left| \int_{S \times T} W(x, y) dx dy \right|$$
$$= \sup_{f,g:[0,1] \to [0,1]} \left| \int_{[0,1]^2} W(x, y) f(x) g(y) dx dy \right|,$$

where the supremum in the last definition is taken over all measurable functions and $S^c = [0, 1] \setminus S$.

The centrality of the cut norm in this theory relies on the relation that it shares with the so-called *cut metric*, that we now introduce. Recall that a map $\phi : [0, 1] \to [0, 1]$ is *measure preserving* if $\phi^{-1}(X)$ is measurable for every measurable $X \subseteq [0, 1]$, and $|\phi^{-1}(X)| = |X|$ (where $|\cdot|$ is the Lebesgue measure on $[0, 1]$). Measure-preserving transformations are a generalization of relabeling of the nodes of a graph. In fact, two different labelings of the nodes of $G$ are related by a permutation, so that two graphons $W_{G,1}, W_{G,2}$ correponding to two labelings of $G$ (cf. (6)) are related by a measure preserving map $\phi$ on $[0, 1]$, such that

$$W_{G,1}^\phi = W_{G,2},$$

where $W^\phi(x, y) = W(\phi(x), \phi(y))$, for every $(x, y) \in [0, 1]^2$ and for $W \in \mathcal{W}$, and $\phi$ is a permutation of the intervals $I_i^n$ in (6).

DEFINITION 2 ([6]). *For every $U, W \in \mathcal{W}$, the* cut distance *between $U$ and $W$, $\delta_\Box(U, W)$, is defined as*

$$\delta_\Box(U, W) = \inf_{\phi, \psi} ||U^\phi - W^\psi||_\Box,$$

*where the infimum is taken over all the measure-preserving maps $\phi, \psi : [0, 1] \to [0, 1]$ (the same holds when the infimum is restricted to bijective measure-preserving maps.)*

The cut distance defined above is actually a pre-metric, because $\delta_\Box(U, W) = 0$ may hold even when $U$ and $W$ are not almost everywhere equal (for further details see for instance [6]). The following result is fundamental.

THEOREM 3 ([27]). *The metric space $(\mathcal{W}_0, \delta_\Box)$ is compact (provided the identification of graphons with $\delta_\Box = 0$).*

Left convergence can be rephrased in terms of graphons, defining, for $W \in \mathcal{W}$, the following graph parameter (i.e., functional defined on graph spaces):

(8) $$t(F, W) = \int_{[0,1]^k} \prod_{ij \in E(F)} W(x_i, x_j) dx_1...dx_k,$$

for any simple graph $F$ with $V(F) = \{1, ..., k\}$ (cf. [6]). This notation is consistent with that of homomorphism densities in (5) since, for every pair of graphs $G$ and $F$,

$$t(F, G) = t(F, W_G),$$

where the left and right hand sides are (5) and (8) respectively, and $W_G$ is the graphon associated to $G$ as in (6).

The following theorem (in the version of [6]) unravels the relation between left convergence and convergence with respect to the cut metric.

THEOREM 4 ([6]). *Let $\{W_n\}_{n \in \mathbb{N}}$ be a sequence of graphons in $\mathcal{W}_0$. Then the following conditions are equivalent:*
- *$t(F, W_n)$ converges as $n \to \infty$, for every finite simple graph $F$;*
- *$\{W_n\}_{n \in \mathbb{N}}$ is a Cauchy sequence in $\delta_\Box$;*
- *there exists $W \in \mathcal{W}_0$ such that $t(F, W_n) \to t(F, W)$ as $n \to \infty$, for every simple graph $F$;*
- *there exists $W \in \mathcal{W}_0$ such that $\delta_\Box(W_n, W) \to 0$.*

Since we can associate a graphon to every simple graph, this implies the following fundamental result.

THEOREM 5 ([6]). *Let $\{G_n\}_{n \in \mathbb{N}}$ be any left convergent sequence of simple graphs. Then there exists a graphon $W \in \mathcal{W}_0$ such that $\delta_\Box(W_{G_n}, W) \to 0$ as $n \to \infty$. The limit is unique up to identification of graphons with cut distance equal to 0.*

*Conversely, every $W \in \mathcal{W}_0$ arises as the limit in the cut metric of a convergent sequence of such graphs.*

Convergence in the cut metric avoids the possible ambiguity that may arise from the non-uniqueness of the association between graphs and graphons, and the limit of any converging graph sequence as in Theorem 5 is always well-defined. However, the cut norm is easier to use in applications, and the following result allows to use it to characterize convergence.

THEOREM 6 ([6]). *Let $\{G_n\}_{n\in\mathbb{N}}$ be a sequence of simple graphs. If*

$$\delta_\square(W_{G_n}, W) \to 0$$

*for some $W \in \mathcal{W}_0$, then there exists a relabelling of the nodes such that the resulting sequence $\{G_{n'}\}_{n'\in\mathbb{N}}$ of labeled graphs converges in the cut norm to $W$; i.e.,*

$$||W_{G_{n'}} - W||_\square \to 0.$$

**Remark 1** (Weak convergence vs convergence in the cut norm). Considering $\mathcal{W}_0$ as a subset of $L^1([0,1]^2)$, the Dunford-Pettis theorem guarantees that $\mathcal{W}_0$ is weakly sequentially compact. However, weak convergence in $L^1([0,1]^2)$ is not appropriate when dealing with graph sequences. Indeed, the homomorphism densities (8) are not continuous in the weak topology, so if $W_n \rightharpoonup W$ in $L^1([0,1]^2)$, this does not imply in general that $t(F, W_n) \to t(F, W)$ for every simple graph $F$. We refer to Appendix F in [25] for further details on this point.

The basic difference between the topology induced by the cut metric and the weak topology in $L^1([0,1]^2)$ is a uniformity requirement. Indeed, $W_n \rightharpoonup W$ if and only if for all fixed $S, T \subseteq [0,1]$, $\int_{S\times T} W_n \to \int_{S\times T} W$. Now, the definition of the cut norm (7) uses the same sets, but now $||W_n - W||_\square \to 0$ if and only if $\int_{S\times T} W_n \to \int_{S\times T} W$ uniformly for all $S, T \subseteq [0,1]$. Clearly, this implies that every sequence of graphons $\{W_n\} \subset \mathcal{W}_0$ converging in the cut norm necessarily also converges in the weak topology of $L^1([0,1]^2)$. Moreover, strong $L^1([0,1]^2)$ convergence is strictly stronger than convergence in the cut norm (see [10] Section 8.3).

From now on we restrict to sequences of graphs with a nearly maximal number of connections, as in the next definition. However, note that the results described above do not require any assumption on the growth rate of the number of edges of $G_n$ as $n \to +\infty$.

DEFINITION 7 (dense graph sequences[2]). *A graph sequence $\{G_n\}_{n\in\mathbb{N}}$ is said to be a* dense graph sequence *when $|V(G_n)| = n$ and*

$$\limsup_{n\to\infty} \frac{n^2}{|E(G_n)|} < +\infty.$$

All graph sequences considered in this paper will be assumed to be dense graph sequences. A trivial example is the sequence of complete graphs $\{K_n\}_{n\in\mathbb{N}}$, with $|V(K_n)| = n$ and $|E(K_n)| = n(n-1)/2$ and all possible connections are present. Indeed, if the edge density $|E(G_n)|/n^2$ tends to 0, the homomorphism densities $t(F, G_n)$ of every finite simple graph $F$ tend to 0 as well, and left convergence does not provide any information on the asymptotic behavior of the sequence.

---

[2] We follow here the custom in the literature on graph sequences and use the term 'dense' to refer to the order of growth of the number of edges of the graphs, rather than to the density of the sequence in the topological sense.

**3. Young measures and weak limits of functions.** We recall some facts about Young measures on $[0,1] \times \mathbb{R}$ (see for instance [28], [29]). There are many equivalent definitions of Young measures, one of which is the following [31]: a *Young measure* $\nu$ on $[0,1] \times \mathbb{R}$ is a family $\{\nu_x\}_{x \in [0,1]} \subset \mathcal{P}(\mathbb{R})$ of probability measures on $\mathbb{R}$, indexed by the elements of $[0,1]$, satisfying the property that

$$(9) \qquad x \mapsto \int_{\mathbb{R}} f(\lambda) d\nu_x(\lambda)$$

is measurable as a function $[0,1] \to \mathbb{R}$, for every continuous and bounded $f \in C^b(\mathbb{R})$. We refer to [30] for some references on applications of Young measures to problems in Continuum Mechanics.

We denote by $\mathcal{Y}([0,1], \mathbb{R})$ the set of Young measures on $[0,1]$. Special cases are Young measures corresponding to functions: if $u$ is a measurable function $u : [0,1] \to \mathbb{R}$, we let

$$(10) \qquad \nu_x = \delta_{u(x)}, \quad \text{for every } x \in [0,1],$$

so that the function in (9) reduces to $x \to f(u(x))$.

We say that a sequence of Young measures $\{\nu^n\}_{n \in \mathbb{N}}$ *converges narrowly* to $\nu \in \mathcal{Y}([0,1], \mathbb{R})$ if, for every fixed continuous and bounded $f \in C^b(\mathbb{R})$, the real function in (9) is weakly-$*$ convergent ([29], [32]):

$$(11) \qquad \int_{\mathbb{R}} f(\lambda) d\nu_x^n(\lambda) \overset{*}{\rightharpoonup} \int_{\mathbb{R}} f(\lambda) d\nu_x(\lambda) \qquad \text{in } L^\infty([0,1]).$$

We will use the following lemma on product Young measures (Theorem 4.17 in [29]).

LEMMA 8. *Let $\{\nu^n\}_{n \in \mathbb{N}}$ be a sequence of Young measures in $\mathcal{Y}([0,1], \mathbb{R})$ converging narrowly to $\nu \in \mathcal{Y}([0,1], \mathbb{R})$. Then the tensor product $\{\nu^n \otimes \nu^n\}_{n \in \mathbb{N}}$ converges narrowly to $\nu \otimes \nu \in \mathcal{Y}([0,1]^2, \mathbb{R}^2)$, i.e.,*

$$\int_{\mathbb{R}^2} f(\lambda, \mu) d\nu_x^n(\lambda) d\nu_y^n(\mu) \overset{*}{\rightharpoonup} \int_{\mathbb{R}^2} f(\lambda, \mu) d\nu_x(\lambda) d\nu_y(\mu) \qquad \text{in } L^\infty([0,1]^2),$$

*for every $f \in C^b([0,1]^2)$.*

The following fundamental result states that boundedness in $L^1$ of sequences of functions implies that the corresponding sequences of Young measures in (10) have subsequences that admit a narrow limit $\nu$, called the *Young measure associated to the sequence* (see for instance [28]).

THEOREM 9 (Prohorov's Theorem). *Let $\{u_n\}_{n \in \mathbb{N}}$ be a norm-bounded sequence in $L^1([0,1], \mathbb{R})$, and denote by $\{\nu^n\}_{n \in \mathbb{N}}$ the corresponding sequence of Young measures as in (10). Then there exists a subsequence $\{u_{n_k}\}_{k \in \mathbb{N}}$ and a Young measure $\nu$ such that $\nu^{n_k} \to \nu$ narrowly as $k \to \infty$.*

**3.1. Functions with finite range and Young measures.** In what follows we shall work with functions taking values in a finite set of numeric labels

$$\mathcal{L} = \{\ell_1, \ldots, \ell_N\} \subset \mathbb{R},$$

and denote by

$$X = \{u : [0,1] \to \mathcal{L} : u \text{ measurable}\} \subset L^\infty([0,1]),$$

the set of such functions on $[0,1]$. Since $X$ is bounded in $L^\infty([0,1])$, it is relatively compact in the weak-$*$ topology of $L^\infty$. Hence, every sequence of such functions admits a weak limit, but since $X$ is not weak-$*$ closed, the limit takes values in the convex hull of $\mathcal{L}$.

Young measures corresponding to functions in $X$ by (10) are such that $\nu_x$ is supported in $\mathcal{L}$ for almost every $x \in [0,1]$. In general, we denote by

$$\mathcal{X} \subset \mathcal{Y}([0,1], \mathbb{R})$$

the set of Young measures such that $\nu_x$ has support in $\mathcal{L}$ for almost every $x \in [0,1]$. Such measures admit the representation

(12) $$\nu \in \mathcal{X} \quad \Leftrightarrow \quad \nu_x = \sum_{k=1}^{N} \theta_k(x) \delta_{\ell_k} \text{ for almost all } x \in [0,1],$$

with $\theta_k : [0,1] \to [0,1]$ measurable functions such that $\sum_{k=1}^{N} \theta_k(x) = 1$. If $\nu$ is the Young measure corresponding to a function in $X$ as in (10), then $\theta_k : [0,1] \to \{0,1\}$ for every $k = 1, \ldots, N$ and for almost every $x \in [0,1]$.

The following result is a consequence of Theorem 4.12 of [32], and ensures that $\mathcal{X}$ is closed as a subset of $\mathcal{Y}([0,1], \mathbb{R})$ with respect to the narrow topology.

PROPOSITION 10. *Let $\{\nu^n\}_{n \in \mathbb{N}} \subset \mathcal{X}$ be a sequence of Young measures converging narrowly to $\nu \in \mathcal{Y}([0,1], \mathbb{R})$. Then $\nu \in \mathcal{X}$.*

**Remark 2.** Since $\mathcal{X}$ is closed with respect to narrow convergence, sequences in $\mathcal{X}$ and their limits are characterized by the measurable functions $\theta_k$ in formula (13). Granted (11) and the identity

$$\int_{\mathbb{R}} f(\lambda) d\nu_x(\lambda) = \sum_{k=1}^{N} f(\ell_k) \theta_k(x),$$

for every $f : \mathcal{L} \to \mathbb{R}$, narrow convergence in $\mathcal{X}$ is equivalent to weak-$*$ convergence in $L^\infty([0,1], \mathbb{R})$ for each $\theta_k$. More precisely, denoting by

$$\Delta = \left\{ \mathbf{y} \in \mathbb{R}^N : \sum_{k=1}^{N} y_k = 1 \text{ and } y_k \geq 0 \ \forall k \right\}$$

the simplex of probability vectors on $\mathcal{L}$, and by

$$\boldsymbol{\theta} = (\theta_1, \ldots, \theta_N) : [0,1] \to \Delta,$$

we have the following equivalence statement between sequences of Young measures in $\mathcal{X}$ and the associated probability vectors:

$$\nu^n \to \nu \text{ in } \mathcal{X} \quad \Leftrightarrow \quad \theta_k^n \overset{*}{\rightharpoonup} \theta_k \text{ in } L^\infty([0,1]), \qquad k = 1, \ldots, N,$$

which is in turn equivalent to $\boldsymbol{\theta}^n \overset{*}{\rightharpoonup} \boldsymbol{\theta}$ in $L^\infty([0,1], \Delta)$.

**Remark 3** (Spin functions and measures). When $\mathcal{L} = \{-1, 1\}$ we shall refer to functions with range in $\mathcal{L}$ as *spin functions*, and to Young measures with support in $\{-1, 1\}$ as *spin Young measures*; such measures admit the representation

(13) $$\nu \in \mathcal{X} \quad \Leftrightarrow \quad \nu_x = \theta(x)\delta_1 + (1 - \theta(x))\delta_{-1} \text{ for almost all } x \in [0,1],$$

with $\theta : [0,1] \to [0,1]$ a measurable function.

If $\nu^n = \{\nu^n_x\}_{x \in [0,1]} \in \mathcal{X}$ is the spin Young measure associated to a spin function $u_n$ as in (10) and $\theta^n : [0,1] \to \{0,1\}$ is as in (13), for every $n \in \mathbb{N}$, we have that

$$u_n(x) = \int_\mathbb{R} \lambda d\nu^n_x(\lambda) = 2\theta^n(x) - 1.$$

Hence, we have the following equivalence statement between sequences of spin functions and the associated spin measures:

$$u_n \overset{*}{\rightharpoonup} u \text{ in } L^\infty([0,1]) \quad \Leftrightarrow \quad \nu^n \to \nu \text{ in } \mathcal{X} \quad \Leftrightarrow \quad \theta^n \overset{*}{\rightharpoonup} \theta \text{ in } L^\infty([0,1]).$$

**3.2. Γ-convergence.** If $U$ is a metrizable topological space, $F_n : U \to [0,+\infty]$, $n \in \mathbb{N}$, a sequence of functionals and $F : U \to [0,+\infty]$, we say that $\{F_n\}_{n\in\mathbb{N}}$ Γ-*converges* to $F$ if for every $u \in U$

(14)
$$\begin{array}{l}(i) \text{ for every } \{u_n\}_{n\in\mathbb{N}} \text{ such that } u_n \to u, \text{ then } F(u) \leq \liminf_n F_n(u_n), \\ (ii) \text{ there exists } \{u_n\}_{n\in\mathbb{N}}, \text{ with } u_n \to u, \text{ such that } F(u) = \limsup_n F_n(u_n).\end{array}$$

The sequence in (ii) is called a *recovery sequence* for $F(u)$.

In what follows we shall assume that either $U$ is a bounded subset of $L^\infty([0,1])$ with the weak-$*$ convergence, or $U = \mathcal{Y}([0,1],\mathbb{R})$ with the narrow convergence. In these cases, thanks to the compactness of the space domains, the Γ-limit of a sequence of functionals has the property that sequences of minimizers of the functionals $F_n$ converge to an absolute minimizer of $F$.

**4. Γ-convergence of the cut functional on dense graph sequences.** Consider a dense graph sequence $\{G_n\}_{n\in\mathbb{N}}$ and, for every $n$, define the discrete cut functional

(15) $$\tilde{F}_n(\tilde{u}) = \frac{1}{n^2} \sum_{i,j \in V(G_n)} A^n_{ij} f(\tilde{u}(i), \tilde{u}(j)),$$

where $A^n_{ij}$ is the adjacency matrix of $G_n$, $\tilde{u}$ is a discrete function $\tilde{u} : V(G_n) \to \mathcal{L}$ and $f$ is a real function defined on the set $\mathcal{L}^2$.

Note that the choice $\mathcal{L} = \{-1,1\}$ and $f(\lambda, \mu) = |\lambda - \mu|^2$ yields the cut functional (1).

Each discrete function $\tilde{u}$ can be extended to a piecewise-constant function in $L^\infty([0,1])$ with values in $\mathcal{L}$ by

(16) $$u : [0,1] \to \mathcal{L}, \quad u(x) = \tilde{u}(i) \quad \text{if } x \in I^n_i,$$

where $I^n_i = (\frac{(i-1)}{n}, \frac{i}{n}]$ is the interval of $[0,1]$ corresponding to the $i$-th node (cf. (6)). The spaces of piecewise-constant functions associated to the subdivision $\{I^n_i\}_{i=1,\ldots,n}$ will be denoted by

$$X_n = \{u : [0,1] \to \mathcal{L} : u \text{ constant on } I^n_i\} \subset X.$$

We denote by $\{W_n\}_{n\in\mathbb{N}}$ the sequence of piecewise-constant graphons associated to $\{G_n\}_{n\in\mathbb{N}}$ by (6), and extend the discrete cut functional (15) to $L^\infty([0,1])$ by

(17) $$F_n(u) = \begin{cases} \int_{[0,1]^2} W_n(x,y) f(u(x), u(y)) dx dy & \text{if } u \in X_n \\ +\infty & \text{if } u \in L^\infty([0,1]) \setminus X_n, \end{cases}$$

so that
$$F_n(u) = \tilde{F}_n(\tilde{u}) \qquad \text{for } u \in X_n,$$
when $u$ and $\tilde{u}$ are related by (16).

We may assume that the sequence $\{G_n\}_{n\in\mathbb{N}}$ is convergent in the sense of Theorem 4, so that $\delta_\square(W_n, W) \to 0$ as $n \to \infty$. By Theorem 6, we can also assume that, modulo relabeling, the sequence $\{W_n\}_{n\in\mathbb{N}}$ converges to $W$ in the cut norm, i.e.,

(18) $$\|W_n - W\|_\square \to 0$$

as $n \to +\infty$. In order to compute the $\Gamma$-limit of the sequence of functionals $\{F_n\}_{n\in\mathbb{N}}$, since $X$ is not closed, we extend it to the set $\mathcal{X}$ of Young measures, which is closed under narrow convergence. Our main result follows from a basic continuity result involving sequences of cut functionals.

LEMMA 11. *Assume that $W_n \to W \in \mathcal{W}_0$ in the cut norm, and let $\{u_n\}_{n\in\mathbb{N}}$ be a sequence of piecewise-constant functions with $u_n \in X_n$. Let $\nu \in \mathcal{X}$ be the Young measure associated to the sequence $\{u_n\}_{n\in\mathbb{N}}$ by Prohorov's theorem. Then, up to a subsequence,*

(19) $$\lim_{n\to+\infty} \int_{[0,1]^2} W_n(x,y) f(u_n(x), u_n(y)) dx dy$$
$$= \int_{[0,1]^2} W(x,y) \Big( \int_{\mathbb{R}^2} f(\lambda, \mu) d\nu_x(\lambda) d\nu_y(\mu) \Big) dx dy$$
$$= \sum_{h,k=1}^{N} f(\ell_h, \ell_k) \int_{[0,1]^2} W(x,y) \theta_h(x) \theta_k(y) dx dy,$$

*where $\nu$ and $\theta_k$ are related by (12).*

*Proof.* We will argue in terms of the weight functions $\theta_k$ in the representation (12). Let $g_n(x,y) := f(u_n(x), u_n(y))$. Denoting by $\nu^n$ the Young measure corresponding to $u_n$ by (10), and by $\theta_k^n$ the associated weight functions, we can write

(20) $$g_n(x,y) = \int_{\mathbb{R}^2} f(\lambda,\mu) d\nu_x^n(\lambda) d\nu_y^n(\mu) = \sum_{h,k=1}^{N} f(\ell_h, \ell_k) \theta_h^n(x) \theta_k^n(y). \qquad \square$$

Now, recalling that $\nu^n \to \nu$ narrowly or, equivalently, $\theta_k^n \overset{*}{\rightharpoonup} \theta_k$ for every $k = 1, \ldots, N$, we have

$$g_n(x,y) \overset{*}{\rightharpoonup} g(x,y) := \int_{\mathbb{R}^2} f(\lambda,\mu) d\nu_x(\lambda) d\nu_y(\mu) = \sum_{h,k=1}^{N} f(\ell_h, \ell_k) \theta_h(x) \theta_k(y)$$

in $L^\infty([0,1]^2)$. This fact follows either from Lemma 8 on the convergence of product measures, or by working directly on the product $\theta_h(x)\theta_k(y)$. Since

$$\left| \int_{[0,1]^2} W_n(x,y) g_n(x,y) dx dy - \int_{[0,1]^2} W(x,y) g(x,y) dx dy \right|$$
$$\leq \left| \int_{[0,1]^2} (W_n(x,y) - W(x,y)) g_n(x,y) dx dy \right|$$
$$+ \left| \int_{[0,1]^2} W(x,y)(g_n(x,y) - g(x,y)) dx dy \right|,$$

recalling that, for every $n$, $\theta_k^n$ takes values in $\{0,1\}$ and letting

$$S_{n,k} = \{x \in [0,1] : u_n(x) = \ell_k\} = \{x \in [0,1] : \theta_k^n(x) = 1\},$$

we have by (7)

$$\left| \int_{[0,1]^2} (W_n(x,y) - W(x,y))g_n(x,y)dxdy \right|$$
$$\leq \sum_{h,k=1}^{N} f(\ell_h, \ell_k) \left| \int_{S_{n,h} \times S_{n,k}} (W_n(x,y) - W(x,y))dxdy \right| \leq C\|W_n - W\|_\square,$$

with $C > 0$ a suitable constant, so that convergence in the cut norm of $\{W_n\}_{n \in \mathbb{N}}$ together with the weak-$*$ convergence of $\{g_n\}_{n \in \mathbb{N}}$ yield the thesis.

**Remark 4.** Uniformity (with respect to the integration domains $S_n$) of the convergence in cut norm plays a fundamental role in the above proof: if $W_n \rightharpoonup W$ only, the result does not hold. As a counterexample, consider the sequence of checkerboard graphons corresponding to a sequence of complete bipartite graphs (cf. Section 5.3.1)

$$W_n(x,y) = \begin{cases} 1 & (x,y) \in (S_n \times S_n^c) \cup (S_n^c \times S_n), \\ 0 & \text{otherwise,} \end{cases}$$

where $S_n = \bigcup_{k=0}^{n-1}[\frac{2k}{2n}, \frac{2k+1}{2n}]$: this sequence converges weakly in $L^1([0,1]^2)$ to the constant graphon $W(x,y) = 1/2$, but does not converge in the cut norm, since

$$\|W_n - 1/2\|_\square \geq \int_{S_n \times S_n^c} (W_n(x,y) - 1/2)dxdy = \frac{1}{2}|S_n \times S_n^c| = \frac{1}{8} > 0.$$

Let $N = 2$, $\mathcal{L} = \{-1, +1\}$, and consider the sequence of spin functions

$$u_n(x) = \begin{cases} 1 & x \in S_n, \\ -1 & x \in S_n^c, \end{cases} \Leftrightarrow \theta^n(x) = \begin{cases} 1 & x \in S_n, \\ 0 & x \in S_n^c, \end{cases}$$

such that $\theta^n \stackrel{*}{\rightharpoonup} 1/2$ in $L^\infty([0,1])$, and choosing $f(1,1) = f(-1,-1) = 0$ and $f(1,-1) = f(-1,1) = 1$, then

$$g_n(x,y) = \theta^n(x)(1 - \theta^n(y)) + \theta^n(y)(1 - \theta^n(x)) = \begin{cases} 1 & (x,y) \in (S_n \times S_n^c) \cup (S_n^c \times S_n), \\ 0 & \text{otherwise,} \end{cases}$$

that weak-$*$ converges to $g(x,y) = 1/2$ in $L^\infty([0,1]^2)$. Then

$$\int_{[0,1]^2} W_n(x,y)g_n(x,y)dxdy = |(S_n \times S_n^c) \cup (S_n^c \times S_n)| = \frac{1}{2}$$
$$\neq \frac{1}{4} = \int_{[0,1]^2} W(x,y)g(x,y)dxdy,$$

so that (19) does not hold.

We define the functionals on Young measures

$$(21) \quad I(\nu) = \begin{cases} \int_{[0,1]^2} W(x,y) \Big( \int_{\mathbb{R}^2} f(\lambda,\mu) d\nu_x(\lambda) d\nu_y(\mu) \Big) dx dy & \nu \in \mathcal{X} \\ +\infty & \nu \in \mathcal{Y}([0,1],\mathbb{R}) \setminus \mathcal{X}, \end{cases}$$

and
(22)
$$I_n(\nu) = \begin{cases} \int_{[0,1]^2} W_n(x,y) \Big( \int_{\mathbb{R}^2} f(\lambda,\mu) d\nu_x(\lambda) d\nu_y(\mu) \Big) dx dy & \nu \in \mathcal{X}_n \\ +\infty & \nu \in \mathcal{Y}([0,1],\mathbb{R}) \setminus \mathcal{X}_n, \end{cases}$$

where $\mathcal{X}_n$ is the set of Young measures corresponding to piecewise-constant functions of finite range in $X_n$, i.e.,

$$\mathcal{X}_n = \{\nu \in \mathcal{X} : \nu_x = \delta_{u(x)}, u \in X_n, \text{ for all } x \in [0,1]\}.$$

Note that

$$I_n(\nu) = F_n(u) \quad \text{when} \quad \nu_x = \delta_{u(x)}.$$

Our main result states that the sequence of cut functionals (22) associated to a dense graph sequence $\Gamma$-converges to $I$, and is a consequence of the continuity result of Lemma 11.

THEOREM 12. *Let $\{G_n\}_{n \in \mathbb{N}}$ be a dense graph sequence, let $W \in \mathcal{W}_0$ its limit graphon, and assume that $\|W_n - W\|_\square \to 0$, where $W_n$ is given by (6). Let $F_n$ be the cut functional (17) associated to $G_n$, for every $n$. Then the sequence $\{F_n\}_{n \in \mathbb{N}}$ $\Gamma$-converges to $I$, in the sense that*

$$(23) \quad \Gamma\text{-}\lim_{n \to \infty} I_n = I,$$

*with respect to narrow convergence in $\mathcal{Y}([0,1],\mathbb{R})$, where $I$ and $I_n$ are given by (21) and (22).*

*Proof.* It is enough to verify that *(i)* and *(ii)* in (14) hold for $\nu \in \mathcal{X}$: in fact, for $\nu \in \mathcal{Y}([0,1],\mathbb{R}) \setminus \mathcal{X}$, we have that eventually $I(\nu) = I_n(\nu^n) = +\infty$ for every $\nu^n \to \nu$, and *(i)* and *(ii)* hold trivially.

Assume now that $\nu^n \to \nu \in \mathcal{X}$: then, if $\{\nu^n\}$ contains a subsequence $\{\nu^{n_k}\}$ such that $\nu^{n_k} \in \mathcal{X}_{n_k}$, by Lemma 11, (21) and (22) we have that $\lim_{k \to \infty} I(\nu^{n_k}) = I(\nu)$. For all other $\nu^n$ in the sequence, $I_n(\nu^n) = +\infty$, and we conclude that *(i)* holds.

To prove *(ii)*, we first show that for every $\nu \in \mathcal{X}$, we can construct a sequence $\{\nu^n\}_{n \in \mathbb{N}}$ such that $\nu^n \in \mathcal{X}_n$ and $\nu^n \to \nu \in \mathcal{X}$. To see this, recall that every Young measure $\nu \in \mathcal{X}$ admits the representation (12), with $\boldsymbol{\theta} \in L^\infty([0,1],\Delta)$. Note that the simplex $\Delta$ is the closed convex hull of its vertices $K = \{\mathbf{e}_1, \ldots, \mathbf{e}_N\}$, with $(\mathbf{e}_i)$ the canonical basis of $\mathbb{R}^N$. Define

$$\begin{aligned}
Z_n =& \{\boldsymbol{\theta} \in L^\infty([0,1],K) : \boldsymbol{\theta} \text{ is constant on } I_i^n\} \\
(24) \quad =& \{\boldsymbol{\theta} \in L^\infty([0,1],\Delta) : \theta_k(x) \in \{0,1\} \text{ for almost all } x \in [0,1] \\
& \qquad \text{and } \theta_k \text{ is constant on all } I_i^n, k = 1, \ldots, N\},
\end{aligned}$$

so that

$$\nu \in \mathcal{X}_n \quad \Leftrightarrow \quad \boldsymbol{\theta} \in Z_n.$$

Given the equivalence between narrow and weak-$*$ convergence, it is enough to prove that for every piecewise constant $\boldsymbol{\theta} \in L^\infty([0,1], \Delta)$ there exists a sequence $\{\boldsymbol{\theta}^n\}_{n \in \mathbb{N}}$ with $\boldsymbol{\theta}^n \in Z_n$ such that $\boldsymbol{\theta}^n \stackrel{*}{\rightharpoonup} \boldsymbol{\theta}$. This is done by a construction analogous to that of the density of Young measures corresponding to (piecewise-constant) functions in the set of all Young measures. In this case the argument is simpler, since we can use the representation in terms of $\boldsymbol{\theta}$. We first note that an approximation argument allows to restrict to $\boldsymbol{\theta}$ piecewise constant and $\boldsymbol{\theta}([0,1]) \subset \mathbb{Q}^N$. Then we may explicitly construct functions $\boldsymbol{\theta}^n$ which give a recovery sequence. For each $\mathbf{t} = (t_1, \ldots, t_N) \in \boldsymbol{\theta}([0,1])$ we denote by $I_{\mathbf{t}}$ the set (union of intervals) where $\boldsymbol{\theta}(x) = \mathbf{t}$. For such $\mathbf{t}$ we may write $t_k = p_k/q$ for some $q$ and $p_k \in \mathbb{N}$, and define

$$\boldsymbol{\theta}^n(x) = \mathbf{e}_k$$

if $x \in I_n^i$, $i/n \in I_{\mathbf{t}}$, and $\sum_{j=1}^{k-1} p_j < i \leq \sum_{j=1}^{k} p_j$ modulo $q$. Here we use the convention that $\sum_{j=1}^{0} p_j = 0$. Note that the condition $\mathbf{t} \in \Delta$ implies that $\sum_{j=1}^{N} p_j = q$, so that the funcion $\boldsymbol{\theta}^n$ is well defined. Moreover, by construction $\boldsymbol{\theta}^n \in Z_n$, and $\boldsymbol{\theta}^n \stackrel{*}{\rightharpoonup} \mathbf{t}$ on $I_{\mathbf{t}}$. The corresponding sequence of Young measures $\nu^n \to \nu$ is such that $\nu^n \in \mathcal{X}_n$, so that it induces a sequence of piecewise-constant functions $\{u_n\}_{n \in \mathbb{N}}$ with $u_n \in X_n$ (given by $u_n(i) = \ell_k$ if and only if $\boldsymbol{\theta}^n(i/n) = \mathbf{e}_k$) and $\nu$ the Young measure associated to the sequence. Applying Lemma 11 to this sequence we obtain that $\lim_{n \to \infty} I(\nu^n) = I(\nu)$, so that, in particular, the lim sup equality ($ii$) holds. □

When the weak-$*$ limit of a sequence of piecewise-constant functions $\{u_n\}_{n \in \mathbb{N}}$ is still a piecewise-constant function $u$, the $\Gamma$-limit can be expressed in terms of $u$. In fact, since $\nu_x = \delta_{u(x)}$ for almost every $x \in [0,1]$, then

$$I(\nu) = \int_{[0,1]^2} W(x,y) f(u(x), u(y)) dx dy.$$

The limit functional (21) is easier to work with when expressed in terms of the functions $\theta_k$ as defined in (12), and Theorem 12 can be restated in terms of $\Gamma$-convergence in $L^\infty([0,1])$ as follows. Recalling (20), define
(25)
$$J(\boldsymbol{\theta}) = \begin{cases} \sum_{h,k=1}^{N} f(\ell_h, \ell_k) \int_{[0,1]^2} W(x,y) \theta_h(x) \theta_k(y) dx dy & \text{if } \boldsymbol{\theta} \in L^\infty([0,1], \Delta) \\ +\infty & \text{if } \boldsymbol{\theta} \in L^\infty([0,1], \mathbb{R}^N) \setminus L^\infty([0,1], \Delta), \end{cases}$$

and
(26)
$$J_n(\boldsymbol{\theta}) = \begin{cases} \sum_{h,k=1}^{N} f(\ell_h, \ell_k) \int_{[0,1]^2} W_n(x,y) \theta_h(x) \theta_k(y) dx dy & \text{if } \boldsymbol{\theta} \in Z_n \\ +\infty & \text{if } \boldsymbol{\theta} \in L^\infty([0,1], \mathbb{R}^N) \setminus Z_n, \end{cases}$$

where $Z_n$ is defined in (24).

We thus have the following reformulation of Theorem (12).

THEOREM 13. *Let $\{G_n\}_{n \in \mathbb{N}}$ a dense graph sequence as in Theorem 12. Then $\{J_n\}_{n \in \mathbb{N}}$ $\Gamma$-converges, with respect to the weak-$*$ convergence in $L^\infty([0,1], \mathbb{R}^N)$, to the functional $J : L^\infty([0,1], \mathbb{R}^N) \to \mathbb{R} \cup \{+\infty\}$ in (25).*

**4.1. Minimal-cut problems.** A mimimal-cut problem for a graph $G_n$ is the problem of finding a partition $\{S_k^n\}_{k=1,\ldots,N}$ of its nodes into $N$ sets with fixed sizes $|S_k^n| = j_k^n$, with $j_k^n \geq 0$ integers and $\sum_{k=1}^N j_k^n = n$, and with minimum number of edges connecting them. As anticipated in the Introduction, this can be restated as the constrained minimum problem for the cut functional

$$\min_{\tilde{u} \in \tilde{E}^n} \tilde{F}_n(\tilde{u}), \tag{27}$$

where $\tilde{F}_n$ is given by (15) and

$$\tilde{E}^n := \{\tilde{u} : V(G_n) \to \mathcal{L} : |\tilde{u}^{-1}(\ell_k)| = j_k^n, \ k = 1, \ldots, N\}.$$

In fact, for a discrete function $\tilde{u} : V(G_n) \to \mathcal{L}$, let $S_k^n = \tilde{u}^{-1}(\ell_k) \subset V(G_n)$, and rewrite the discrete functional (15) as

$$F_n(\tilde{u}) = \frac{1}{n^2} \sum_{h,k=1}^N \sum_{i \in S_h^n} \sum_{j \in S_k^n} A_{ij} f(\ell_h, \ell_k) = \frac{1}{n^2} \sum_{h,k=1}^N e_{G_n}(S_h^n, S_k^n) f(\ell_h, \ell_k),$$

with $e_{G_n}(S_h^n, S_k^n)$ the number of edges between $S_h^n$ and $S_k^n$. Assuming that $f(\ell_k, \ell_k) = 0$ for $k = 1, \ldots, N$, so that edges connecting nodes with the same label are not penalized, the minimization of the cut functionals is indeed equivalent to minimizing the sizes of the cuts between the partitions.

In terms of piecewise-constant functions, introducing the constraint sets

$$E^n = \left\{u \in L^\infty([0,1]) : |u^{-1}(\ell_k)| = j_k^n/n, \ k = 1, \ldots, N\right\},$$

since (27) has always a solution, also the minimum problem

$$\min_{u \in E^n} F_n(u), \tag{28}$$

has always a solution as well. The minimizers are the piecewise-constant functions in $X_n \cap E^n$ related by (16) to the discrete functions that minimize $\tilde{F}_n(\tilde{u})$. Assume now that $j_k^n/n \to j_k$, and define

$$E = \left\{u \in L^\infty([0,1]) : |u^{-1}(\ell_k)| = j_k, \ k = 1, \ldots, N\right\},$$

and

$$\mathcal{E} = \left\{\nu \in \mathcal{Y}([0,1], \mathbb{R}) : \left|\left\{x \in [0,1] : \int_\mathbb{R} \lambda d\nu_x(\lambda) = \ell_k\right\}\right| = j_k, \ k = 1, \ldots, N\right\},$$

and, finally,

$$H^n = \left\{\boldsymbol{\theta} \in L^\infty([0,1], \Delta) : \int_0^1 \theta_k(x) dx = j_k^n/n, \ k = 1, \ldots, N\right\},$$

$$H = \left\{\boldsymbol{\theta} \in L^\infty([0,1], \Delta) : \int_0^1 \theta_k(x) dx = j_k, \ k = 1, \ldots, N\right\}.$$

The following result provides the basis for the applications presented in the following section.

PROPOSITION 14. *Let $\{G_n\}_{n\in\mathbb{N}}$ be a dense graph sequence and assume that the hypotheses of Theorem 12 are satisfied. Let $j_k^n/n \to j_k$ and $\{u_n\}_{n\in\mathbb{N}}$ be a sequence of minimizers of $F_n$ on $E^n$, i.e., solutions of (28), and $\nu \in \mathcal{Y}([0,1], \mathbb{R})$ the Young measure associated to the sequence. Then $\nu$ is an absolute minimizer of the functional $I$ on $\mathcal{E}$, and there exists the limit*

$$\lim_n \min_{u \in E^n} F_n(u) = \lim_n F_n(u_n) = I(\nu) = \min_{\nu' \in \mathcal{E}} I(\nu'), \tag{29}$$

*with $I$ given by (21). Equivalently,*

$$\lim_n \min_{u \in E^n} F_n(u) = \lim_n F_n(u_n) = J(\boldsymbol{\theta}) = \min_{\boldsymbol{\theta}' \in H} J(\boldsymbol{\theta}'), \tag{30}$$

*with $J$ and $\boldsymbol{\theta}$ defined by (25) and (12).*

*Proof.* We only prove (30). We note that in order to deduce the convergence of minimum problems it suffices to prove the Γ-convergence of the functionals

$$\widetilde{J}_n(\boldsymbol{\theta}) = \begin{cases} J_n(\boldsymbol{\theta}) & \text{if } \boldsymbol{\theta} \in Z_n \cap H^n \\ +\infty & \text{otherwise}, \end{cases} \tag{31}$$

to

$$\widetilde{J}(\boldsymbol{\theta}) = \begin{cases} J(\boldsymbol{\theta}) & \text{if } \boldsymbol{\theta} \in H \\ +\infty & \text{otherwise}. \end{cases} \tag{32}$$

Note that if $\boldsymbol{\theta}^n \in Z_n \cap H^n$ and $\boldsymbol{\theta}^n \stackrel{*}{\rightharpoonup} \boldsymbol{\theta}$, then $\boldsymbol{\theta} \in H$, so that the liminf inequality $(i)$ in (14) is deduced from the corresponding inequality for $J_n$.

In order to prove $(ii)$ in (14), it suffices to show that for every $\widetilde{\boldsymbol{\theta}} \in L^\infty([0,1], \Delta) \cap H$ there exists a sequence $\{\widetilde{\boldsymbol{\theta}}^n\}_{n\in\mathbb{N}}$, with $\widetilde{\boldsymbol{\theta}}^n \in Z_n \cap H^n$, such that $\widetilde{\boldsymbol{\theta}}^n \stackrel{*}{\rightharpoonup} \widetilde{\boldsymbol{\theta}}$, so that $\widetilde{J}_n(\widetilde{\boldsymbol{\theta}}^n) \to \widetilde{J}(\widetilde{\boldsymbol{\theta}})$ by Lemma 11. This is done by modifying a sequence as constructed in the proof of Theorem 12. Indeed, since $\boldsymbol{\theta} \in H$ the sequence $\{\boldsymbol{\theta}^n\}_{n\in\mathbb{N}}$ exhibited there satisfies,

$$\int_0^1 \theta_k^n(x)\,dx \to j_k \quad \text{so that} \quad \int_0^1 \theta_k^n(x)\,dx - \frac{j_k^n}{n} \to 0.$$

We can therefore modify the definition of $\{\boldsymbol{\theta}^n\}_{n\in\mathbb{N}}$ in a set of indices $I_n$ with $|I_n| = o(n)$ in such a way that the resulting $\widetilde{\boldsymbol{\theta}}^n \in H^n$. Since $|\{x \in [0,1] : \widetilde{\boldsymbol{\theta}}^n \neq \boldsymbol{\theta}^n\}| = o(1)$, we still can apply Lemma 11 and obtain the desired equality.

We can now apply the fundamental theorem of the convergence of minima for Γ-converging sequences and the identification of $F_n$ with $J_n$. Namely, consider a sequence of minimizers $u_n$ of $F_n$ on $E^n$, let $\nu^n$ and $\boldsymbol{\theta}^n$ be the corresponding Young measure and mass vector, respectively, and $\nu$ and $\boldsymbol{\theta}$ the Young measure and mass vector associated to the sequence, so that $\nu^n \to \nu$ narrowly and $\boldsymbol{\theta}^n \stackrel{*}{\rightharpoonup} \boldsymbol{\theta}$. Note that $\boldsymbol{\theta}^n$ is a minimizer of $J_n$ on $H^n$. Then by inequality $(i)$ in (14) we have

$$J(\boldsymbol{\theta}) \leq \liminf_{n \to +\infty} J_n(\boldsymbol{\theta}^n) = \liminf_{n \to +\infty} \min_{u \in E^n} F_n.$$

Conversely, since $J(\boldsymbol{\theta}) = \widetilde{J}(\boldsymbol{\theta})$ then by $(ii)$ in (14) there exists $\widetilde{\boldsymbol{\theta}}^n \stackrel{*}{\rightharpoonup} \boldsymbol{\theta}$, so that $\widetilde{J}_n(\widetilde{\boldsymbol{\theta}}^n) \to J(\boldsymbol{\theta})$ and

$$J(\boldsymbol{\theta}) = \lim_{n \to +\infty} \widetilde{J}_n(\widetilde{\boldsymbol{\theta}}^n) \geq \limsup_n \min \widetilde{J}_n = \limsup_{n \to +\infty} \min_{u \in E^n} F_n.$$

Using the two inequalities above we obtain (30). Claim (29) follows in the same way. □

**4.1.1. The graph bisection problem.** The bisection problem for a graph $G_n$, with $n \in 2\mathbb{N}$ even, is the problem of finding a partition $\{S_n, S_n^c\}$ of its nodes with $|S_n| = |S_n^c|$ and with minimum number of edges connecting them. We can choose $\mathcal{L} = \{-1, 1\}$ and $f(\lambda, \mu) = |\lambda - \mu|^2$, so that

$$\int_{\mathbb{R}^2} f(\lambda, \mu) d\nu_x(\lambda) d\nu_y(\mu) = 4\theta(x)(1 - \theta(y)) + 4(1 - \theta(x))\theta(y),$$

and using the symmetry of $W$, we have
(33)
$$J(\theta) = \begin{cases} 8 \int_{[0,1]^2} W(x,y)\theta(x)(1-\theta(y))dxdy & \theta \in L^\infty([0,1],[0,1]) \\ +\infty & \theta \in L^\infty([0,1]) \setminus L^\infty([0,1],[0,1]). \end{cases}$$

Proceeding as above, and defining

$$E = \left\{ u \in L^\infty([0,1]) : \int_0^1 u(x)dx = 0 \right\}, \quad H = \left\{ \theta \in L^\infty([0,1]) : \int_0^1 \theta(x) = \frac{1}{2} \right\},$$

it follows from Proposition 14 that, for a sequence $\{u_n\}$ of minimizers of the cut functionals $F_n$,

(34) $$\min_{u \in E} F_n(u) = F_n(u_n) \to J(\theta) = \min_{\theta' \in H} J(\theta').$$

**5. Examples.** In this section we discuss the properties of the minimizers of the cut functional (33) in some simple cases, and we assume that $n$ is even unless otherwise specified in order that the bisection problem be well posed.

**5.1. The complete graph.** The complete graph $K_n$ on $n$ nodes is the graph that has all possible $n(n+1)/2$ edges connecting all its nodes to each other. A pictorial representation is given in Fig. 1(a).

The limit graphon is the constant graphon $W(x,y) \equiv 1$. Hence, the limit cut functional reduces to

(35) $$J(\theta) = 8 \int_{[0,1]^2} \theta(x)(1-\theta(y))dxdy = 8 \int_0^1 \theta(x)dx \left(1 - \int_0^1 \theta(x)dx\right),$$

for $\theta \in Z$. On the functions $\theta \in H$ satisfying the mass constraint, the functional above is identically constant with $J(\theta) \equiv 2$, so that every such spin measure is a minimizer. When $\theta(x) \in \{0, 1\}$ for almost every $x \in [0, 1]$, the minimizer is a spin function, and the corresponding cut divides the graph in two well-defined regions, where either $u = 1$ or $u = -1$. However, if $\theta(x) \notin \{0, 1\}$ on a set of positive measure, there the cut is 'diffuse', and $\theta(x)$ is the probability that the spin $u(x) = +1$.

As an example, consider the bisection of $K_n$ given by the spin function

(36) $$u_n(x) = \text{sign}(\sin(n\pi x)) \quad x \in [0, 1];$$

the sequence weakly converges to zero, which is not a spin function, and the associated spin Young measure is a mixture of 1 and $-1$ with uniform probability $\theta(x) \equiv 1/2$.

This is consistent with the fact that, in a complete graph, all bisections have the same cut size.

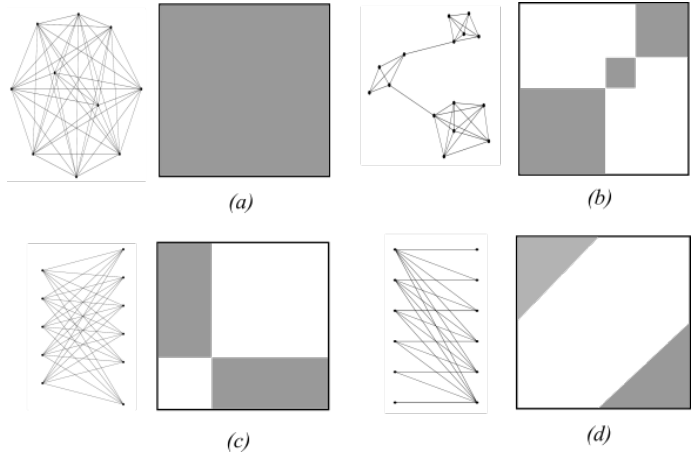

Fig. 1. *Examples of graphs and the support of limit graphons on the square $[0,1]^2$. Darker areas are regions in which $W = 1$. (a): the complete graph; (b) weakly connected complete subgraphs; (c): a complete bipartite graph; (d): the half graph.*

**5.2. Weakly-connected complete subgraphs.** As a generalization of complete graphs, choose strictly positive real numbers $\lambda_1, \ldots, \lambda_N$ such that $\sum_{h=1}^{N} \lambda_h = 1$, let $\lambda_0 = 0$, and define a partition of $\{1, \ldots, n\}$ into $N$ subsets

$$C_k^n = \left\{ \left\lfloor n\sum_{h=0}^{k-1} \lambda_h \right\rfloor + 1, \left\lfloor n\sum_{h=0}^{k-1} \lambda_h \right\rfloor + 2, \ldots, \left\lfloor n\sum_{h=0}^{k} \lambda_h \right\rfloor \right\}, \qquad k = 1, \ldots, N,$$

with $\lfloor x \rfloor$ the largest integer less than or equal to $x$. All pairs of nodes in the same set $C_k^n$ are assumed to be connected, exactly one node in each $C_k^n$ is connected with a node in $C_{k+1}^n$, and no other edge exists. Hence, the sets $C_k^n$ are the sets of nodes of complete subgraphs, each connected by a single edge with its neighbors (see the picture on the left-hand side in Fig. 1(b)). The limit graphon associated to this sequence is

$$W(x,y) = \begin{cases} 1 & \text{if } (x,y) \in C_k \times C_k, \quad k = 1, \ldots, N, \\ 0 & \text{otherwise}, \end{cases}$$

where $C_k = [\sum_{h=0}^{k-1} \lambda_h, \sum_{h=0}^{k} \lambda_h]$, represented in the picture on the right-hand side in Fig. 1(b). The cut functional becomes

$$(37) \qquad J(\theta) = 8 \sum_{h=1}^{k} \int_{C_k} \int_{C_k} \theta(x)(1-\theta(y)) dx dy = 8 \sum_{h=1}^{k} A_k(\lambda_k - A_k),$$

with $A_k = \int_{C_k} \theta(x) dx$, so that $\sum_{k=1}^{N} A_k = 1/2$.

We claim that the minimum-cut problem admits a solution corresponding to a subdivision into regions with vanishing cut size, i.e., a spin function minimizer $\theta(x) \in \{0,1\}$ for all $x \in [0,1]$ such that $J(\theta) = 0$, if and only if there exist $\{k_1, \ldots, k_L\} \subset \{1, \ldots, N\}$ such that

$$(38) \qquad \sum_{j=1}^{L} \lambda_{k_j} = \frac{1}{2}.$$

In fact, if this is the case, the minimizers have the form

$$\theta(x) = \begin{cases} 1 & \text{if } x \in \bigcup_{j=1}^{L} C_{k_j} \\ 0 & \text{otherwise,} \end{cases}$$

or its symmetric $1 - \theta$. The fact that $\theta$ above is a minimizer follows from the fact that $J(\theta) = 0$ and $J$ is non negative. Conversely, assume that $J(\theta) = 0$; then by (37) either $A_k = 0$ or $A_k = \lambda_k$, i.e., either $\theta(x) = 0$ or $\theta(x) = 1$ on $C_k$. Using the mass constraint, we see that the measure of the union of the intervals on which $\theta(x) = 1$ must be $1/2$.

Condition (38) is satisfied for instance if $N$ is even and $\lambda_i = \lambda_j$ for all $i, j$, while it cannot be satisfied if $N$ is odd and $\lambda_i = \lambda_j$ for all $i, j$ as well. In this last case, the minimizer must be such that $\theta(x) \in (0, 1)$ on some interval.

Hence, most subdivisions involve an interface, and there are in general infinitely many solutions, that involve both spin functions (i.e., a sharp subdivision between the partitions) and spin measures. In fact, the minimization of (37) only gives information on $A_i$, and this does not determine $\theta$ uniquely.

As an example, consider a 'dumbbell graphon' with $N = 3$ and $\lambda_3 < \lambda_2 < \lambda_1 < 1/2$ as in Fig. 1(b), so that (37) becomes

$$J(\theta) = A_1(\lambda_1 - A_1) + A_2(\lambda_2 - A_2) + A_3(\lambda_3 - A_3),$$

with the mass constraint $\sum_i A_i = 1/2$, and where

$$0 \leq A_1 \leq \lambda_1, \quad 0 \leq A_2 \leq \lambda_2, \quad \frac{1}{2} - \lambda_3 \leq A_1 + A_2 \leq \frac{1}{2}.$$

The absolute minimizers of the cut functional can be computed by minimizing the function

$$g(A_1, A_2) = A_1(\lambda_1 - A_1) + A_2(\lambda_2 - A_2) + (1/2 - A_1 - A_2)(1/2 - \lambda_1 - \lambda_2 + A_1 + A_2),$$

on the polygon

$$\tilde{C} = \left\{ (A_1, A_2) : 0 \leq A_i \leq \lambda_i, \frac{1}{2} - \lambda_3 \leq A_1 + A_2 \leq \frac{1}{2} \right\}.$$

Since the function $g$ is concave, its minima are reached at some of the vertices of the polygon $\tilde{C}$, i.e.,

$$A = \left(\frac{1}{2} - \lambda_3, 0\right), \quad B = (\lambda_1, 0), \quad C = \left(\lambda_1, \frac{1}{2} - \lambda_1\right),$$
$$D = \left(\frac{1}{2} - \lambda_2, \lambda_2\right), \quad E = (0, \lambda_2), \quad F = \left(0, \frac{1}{2} - \lambda_3\right).$$

We find that $g(A) = g(D)$, $g(B) = g(E)$, $g(C) = g(F)$, and

$$g(A) - g(B) = \left(\frac{1}{2} - \lambda_2\right)(\lambda_1 - \lambda_3), \quad g(A) - g(F) = \left(\frac{1}{2} - \lambda_3\right)(\lambda_1 - \lambda_2),$$
$$g(C) - g(B) = \left(\frac{1}{2} - \lambda_1\right)(\lambda_2 - \lambda_3)$$

so that $g(A) = g(D) > g(F) = g(C) > g(E) = g(B)$, and therefore the absolute minima of $f$ subject to the given constraint are attained at $A_1 = 0$, $A_2 = \lambda_2$, $A_3 = $ $1/2 - \lambda_2$ and $A_1 = \lambda_1$, $A_2 = 0$, $A_3 = 1/2 - \lambda_1$, corresponding to $\theta = 0$ on $C_1$, $\theta = 1$ on $C_2$ and $\theta$ on $C_3$ such that

$$\int_{C_3} \theta(x) dx = \frac{1}{2} - \lambda_2,$$

together with its symmetric $1 - \theta$. The structure of this minimizer is sketched in Figure 2: the interface is restricted to the smallest complete subgraph with nodes in $C_3$, but the interface itself can be either sharp or diffuse.

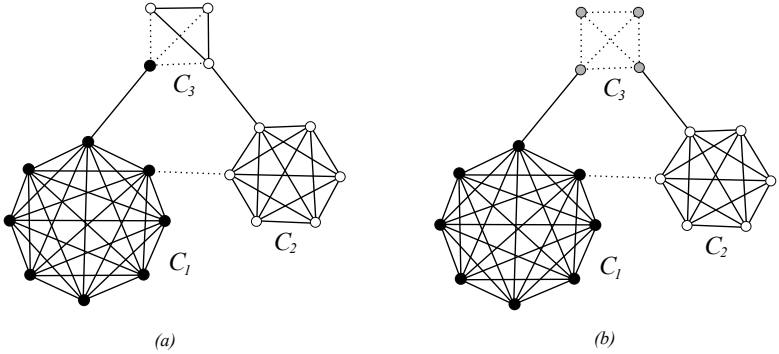

FIG. 2. *The structure of the optimal bipartition in a graph made of three weakly connected complete subgraphs. (a): sharp interface solution; (b): blurry interface solution. Dashed lines are the cuts.*

**5.3. Bipartite graphons.** A graphon $W$ is *bipartite* if there exists a partition $C_1, C_2$ of $[0, 1]$ such that $\operatorname{supp}(W) \subset (C_1 \times C_2) \cup (C_2 \times C_1)$, where $\operatorname{supp}(W)$ is the essential support of $W$. For a bipartite graphon, the cut functional (25) becomes

$$(39) \quad J(\theta) = 8 \int_{C_2} \int_{C_1} W(x,y)[\theta(x)(1 - \theta(y)) + \theta(y)(1 - \theta(x))] dy dx.$$

When $\theta$ corresponds to a spin function, i.e., $\theta(x) \in \{0, 1\}$ almost everywhere, the cut functional has a simple graphical interpretation: let $S_0 \subset [0, 1]$ and $S_1 \subset [0, 1]$ be the sets on which $\theta = 0$ and $\theta = 1$, respectively, and note that the integrand in (39) is non zero only in the checkerboard region $T = (S_0 \times S_1) \cup (S_1 \times S_0)$, so that for spin functions the value of the cut functional is proportional to the area of the region $T \cap \operatorname{supp}(W)$.

**5.3.1. The complete bipartite graph.** We denote by $K_{p,q}$, with $p + q = n$, the bipartite graph with partitions having $p$ and $q$ vertices, such that every node of a partition is connected to all nodes of the other partition, but with no node of its own partition (see the picture on the left-hand side in Fig. 1(c)). Consider a sequence $p_n, q_n$ such that $p_n/n \to \gamma \in [0, 1]$ as $n \to +\infty$, so that the proportions of the nodes in the network that belong to each group are asymptotically constant. Formally, this corresponds to the partition $C_1^n = \{1, \ldots, \lfloor n\gamma \rfloor\}$ and $C_2^n = \{\lfloor n\gamma \rfloor + 1, \ldots, n\}$, and $p_n = \lfloor n\gamma \rfloor$, $q_n = n - p_n$.

The graph sequence $\{K_{p_n, q_n}\}_{n \in 2\mathbb{N}}$ is a dense graph sequence, since the number of edges is $p_n q_n$, which is of order $n^2$. The limit graphon is

$$W(x, y) = \begin{cases} 1 & (x, y) \in ([0, \gamma] \times [\gamma, 1]) \cup ([\gamma, 1] \times [0, \gamma]) \\ 0 & \text{otherwise}, \end{cases}$$

and is represented in the picture on the right-hand side in Fig. 1(c). Note that here $C_1 = [0, \gamma]$ and $C_2 = [\gamma, 1]$. The limit cut functional is

$$J(\theta) = 8 \int_\gamma^1 \int_0^\gamma (\theta(x)(1-\theta(y)) + \theta(y)(1-\theta(x)))dxdy.$$

Letting $A = \int_0^\gamma \theta(x)dx = 1/2 - \int_\gamma^1 \theta(x)dx$, we have that, for every $\theta \in H$,

$$J(\theta) = 8\left[\left(\frac{1}{2} - A\right)(\gamma - A) + A\left(1 - \gamma - \frac{1}{2} + A\right)\right] = 8\left(2A^2 - 2\gamma A + \frac{\gamma}{2}\right),$$

which is minimized by $A = \gamma/2$, i.e., minima are attained whenever half of the nodes in each partition have spin equal to $-1$ and half to $+1$.

Similar to the complete graph, there are infinite minimizers, some corresponding to spin functions, and therefore to actual subdivisions, other which involve 'diffuse' interfaces.

The same result is true also when the two partitions are disconnected, i.e.,

$$W(x,y) = \begin{cases} 1 & \text{if } (x,y) \in \bigcup_{i,j}(C_i \times C'_j) \text{ or } (x,y) \in \bigcup_{i,j}(C'_j \times C_i) \\ 0 & \text{otherwise} \end{cases}$$

where $\{C_i\}_i \cup \{C'_j\}_j$ is a partition of $[0,1]$ such that $|\bigcup_i C_i| = \gamma$.

**5.4. The half-graph.** The half graph $H_{n/2,n/2}$ is a bipartite graph on $n$ nodes with $C_1^n = \{1, \ldots, n/2\}$ and $C_2^n = \{n/2 + 1, \ldots, n\}$ (recall that $n$ is even), such that node $i \in C_1^n$ is connected with node $j \in C_2^n$ if and only if $i \leq j - n/2$ (see the picture on the left-hand side in Fig. 1(d)).

The sequence $\{H_{n/2,n/2}\}$ converges in the cut metric to the graphon $W$ given by

(40) $$W(x,y) = \begin{cases} 1 & \text{if } y + \frac{1}{2} \leq x \text{ or } x + \frac{1}{2} \leq y \\ 0 & \text{otherwise,} \end{cases}$$

represented in the picture on the right-hand side in Fig. 1(d) (cf. [4]).

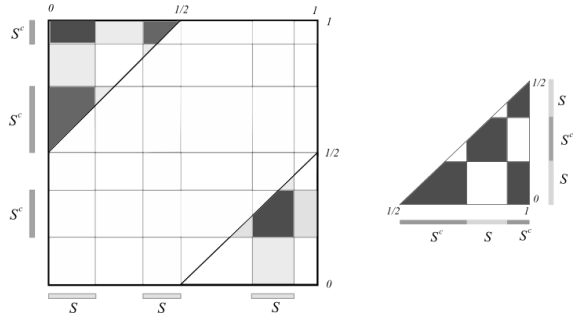

FIG. 3. *Reduction by symmetry for the half graph. $S$ is the set on which $\theta(x) = 1$. The dark areas in the picture on the left correspond to the integral in* (25), *while dark areas in the picture on the right correspond to the integral in* (41).

The $\Gamma$-limit of the cut functional becomes (cf. (39))

(41) $$J(\theta) = 8 \int_{\frac{1}{2}}^1 \int_0^{x-\frac{1}{2}} [\theta(x)(1-\theta(y)) + \theta(y)(1-\theta(x))]\, dy dx.$$

For spin functions $\theta \in \{0, 1\}$, the cut functional above admits a geometric interpretation as the area of the intersection of the checkerboard pattern in Fig. 3(right) with the south-east triangle of the support of the graphon.

**5.4.1. Euler equations for the half-graph.** A necessary condition for a measurable function $\theta \in Z$ to be a minimum of the cut functional (33), subject to the mass constraint, is that the Gateaux differential of the augmented functional $\tilde{J}$ (defined below) is non negative, i.e.,

(42) $$\delta \tilde{J}(\theta, \lambda; \delta\theta) \geq 0,$$

for every $\delta\theta \in L^\infty([0, 1])$ such that

$$\delta\theta(x) = \begin{cases} \geq 0 & \text{if } \theta(x) = 0 \\ \leq 0 & \text{if } \theta(x) = 1 \\ \text{unrestricted} & \text{if } \theta(x) \in (0, 1), \end{cases}$$

where

$$\tilde{J}(\theta, \lambda) = J(\theta) + \lambda \left( \int_0^1 \theta(x) dx - 1/2 \right),$$

and $\lambda$ is the Lagrange multiplier associated to the mass constraint. The Gateaux differential of $\tilde{J}$ is

$$\delta \tilde{J}(\theta, \lambda; \delta\theta) = 8 \int_0^1 \left[ \lambda - \int_0^1 W(x, y)(1 - 2\theta(y)) dy \right] \delta\theta(x) dx.$$

Assume now that there exists $T \subset [0, 1]$, with $|T| > 0$, such that $\theta(x) \in (0, 1)$ for $x \in T$. Then by (42) a necessary condition for $\theta$ to be a minimum is that

(43) $$\int_0^1 W(x, y)(1 - 2\theta(y)) dy = \lambda = \text{constant} \qquad x \in T.$$

Note that (43) is always satisfied for $T = [0, 1]$ and $\theta \equiv 1/2$. However, as shown below, this is not always a minimum of the cut functional, but only a stationary point.

**5.4.2. Minima are spin functions.** We now show that the minima of the cut functional are spin functions. Let $\theta$ be a minimum and assume for contradiction that there exists an interval $R \subset [1/2, 1]$ on which $\theta(x) \in (0, 1)$ for almost every $x \in R$. Then, by (43),

$$\int_0^{x-1/2} (1 - 2\theta(y)) dy = \lambda = \text{constant}, \qquad x \in R,$$

which implies that $1 - 2\theta(x) = 0$ on $R - 1/2$, and in particular that $\theta(x) \in (0, 1)$ for almost every $x \in R - 1/2$. Repeating the argument on $R - 1/2$ it follows that also $1 - 2\theta(x) = 0$ on $R$. Hence, if there exists a minimizer which is not a spin function, it must be identically $\theta(x) = 1/2$ on $R \cup (R - 1/2)$. We claim that there exists a spin function on $R \cup (R - 1/2)$ whose cut size is lower than this. Let $R = [a, b]$, $R - 1/2 = [\hat{a}, \hat{b}]$, $l = b - a$, and for $\lambda \in [0, 1]$ define

$$\tilde{\theta}(x) = \begin{cases} 0 & \text{if } x \in [\hat{a}, \hat{a} + (1-\lambda)l] \cup [a, a + \lambda l] \\ 1 & \text{if } x \in [\hat{a} + (1-\lambda)l, \hat{b}] \cup [a + \lambda l, b] \\ \theta(x) & \text{elsewhere}, \end{cases}$$

and $\bar\theta$ defined by switching 0 and 1 in the above expression.

The values of the cut functional $J(\theta)$ and $J(\tilde\theta)$ are different only in the regions $a, b, c$ in Fig. 4. As to region $a$ we have

$$\int_a^b \int_{\hat a}^{x-1/2} \theta(x)dx = 2l^2, \qquad \int_a^b \int_{\hat a}^{x-1/2} \tilde\theta(x)dx = 8l^2(3\lambda^2 - 4\lambda + 3/2).$$

The minimum of the above expression is attained at $\lambda = 2/3$ and is $4l^2/3$. As to region $b$, we have

$$\int_a^b \int_0^{\hat a} \theta(x)dx = 8l(\mu_1/2 + \mu_2/2), \qquad \int_a^b \int_0^{\hat a} \tilde\theta(x)dx = 8l(\lambda\mu_1 + (1-\lambda)\mu_2),$$

and in region $c$, we have

$$\int_b^1 \int_{\hat a}^{\hat b} \theta(x)dx = 8l(\nu_1/2 + \nu_2/2), \qquad \int_b^1 \int_{\hat a}^{\hat b} \tilde\theta(x)dx = 8l((1-\lambda)\nu_1 + \lambda\nu_2),$$

where $\mu_1 = |\{x \in [0, \hat a] : \theta(x) = 1\}|$, $\mu_2 = |\{x \in [0, \hat a] : \theta(x) = 0\}|$, $\nu_1 = |\{x \in [b, 1] : \theta(x) = 1\}|$, $\nu_2 = |\{x \in [b, 1] : \theta(x) = 0\}|$. Hence, the integrals of $\theta$ and $\tilde\theta$ in the union of regions $b$ and $c$ become respectively

$$8l((\mu_1 + \nu_2)/2 + (\mu_2 + \nu_1)/2), \qquad 8l(\lambda(\mu_1 + \nu_2) + (1-\lambda)(\mu_2 + \nu_1)).$$

On the other hand, we also have that the integral of $\bar\theta$ in the union of regions $b$ and $c$ is

$$8l(\lambda(\mu_2 + \nu_1) + (1-\lambda)(\mu_1 + \nu_2)).$$

Hence, taking $\lambda = 2/3$, either the cut functional involving $\tilde\theta$ or the cut functional involving $\bar\theta$ is smaller than that involving $\theta$, which implies that $\theta$ cannot be a minimizer.

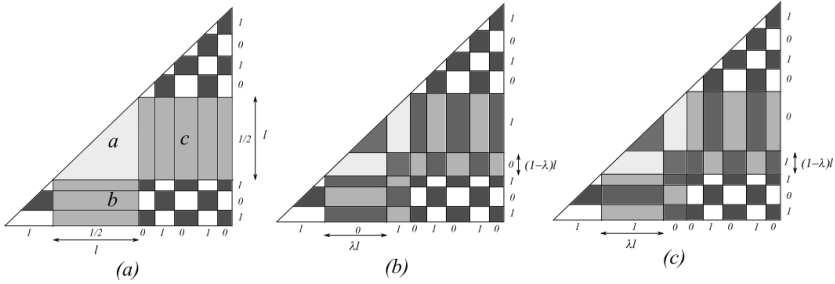

Fig. 4. *Sketch for the proof that the minimum cut for the half graph is a spin function.*

An interesting fact regarding the half graph (and bipartite graphs in general), is that $\theta|_{[0,1/2]}$ and $\theta|_{[1/2,1]}$ are independent, apart from the mass constraint. Integrating (41) by parts, and letting

$$w_1(x) = \int_0^x \theta(\xi)d\xi, \qquad w_2(x) = \int_0^x \theta(\xi + \tfrac{1}{2})d\xi, \qquad \mathbf{w}(x) = (w_1(x), w_2(x)),$$

the cut functional becomes

$$J(\mathbf{w}) = 8 \int_0^{\frac{1}{2}} ((\tfrac{1}{2} - x)w_1'(x) + xw_2'(x) - 2w_2'(x)w_1(x))dx,$$

where $\mathbf{w} = (w_1, w_2) \in W^{1,\infty}([0, \tfrac{1}{2}], \mathbb{R}^2)$ satisfy the boundary conditions

$$w_1(0) = w_2(0) = 0 \qquad w_1(\tfrac{1}{2}) + w_2(\tfrac{1}{2}) = \tfrac{1}{2},$$

and the constraint

$$w_1'(x), w_2'(x) \in [0, 1], \qquad \text{for } x \in [0, 1/2].$$

**5.4.3. The optimal bipartition of the half graph.** The minimizer of the cut functional can be computed by any numerical optimization procedure, and is

$$u(x) = \begin{cases} 1 & x \in [0, 1/6) \cup [1/2, 5/6), \\ -1 & x \in [1/6, 1/2) \cup [5/6, 1], \end{cases}$$

and its symmetric, and the corresponding partition is

$$S = \left[0, \tfrac{1}{6}\right) \cup \left[\tfrac{1}{2}, \tfrac{5}{6}\right), \qquad S^c = \left[\tfrac{1}{6}, \tfrac{1}{2}\right) \cup \left[\tfrac{5}{6}, 1\right].$$

The structure of this minimizer is sketched in Fig. 5: the minimal cut involves a bipartite-complete connection and a half-graph connection between the communities.

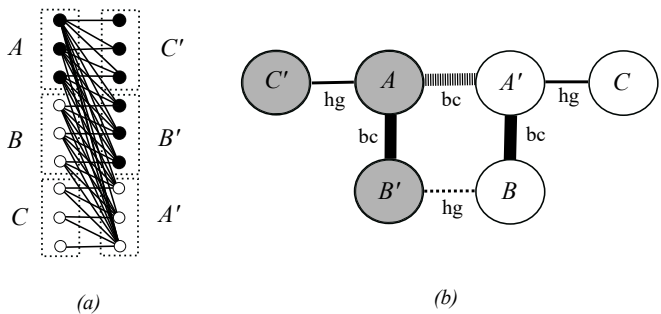

FIG. 5. *The structure of the optimal bipartition in the half graph. (a): standard representation; (b): a schematic representation of the communities: 'hg' stands for a half-graph connectivity between two communities of size k (involving $k(k+1)/2$ edges), while 'bc' stands for a complete bipartite connectivity ($k^2$ edges). Dashed lines are the cuts.*